\documentclass[authoryear,3p, times,10pt,final]{elsarticle}%gulliver

\usepackage{amsmath}
\usepackage{hyperref}

\begin{document}

\begin{frontmatter}

\newtheorem{thm}{Theorem}
\newtheorem{lem}[thm]{Lemma}
\newdefinition{rmk}{Remark}
\newdefinition{note}{Note}
\newproof{pf}{Proof}
\newproof{pot}{Proof of Theorem \ref{thm2}}

%%\usepackage{draftwatermark}
%\usepackage{draftwatermark}
%\SetWatermarkText{Confidential}
%%\SetWatermarkColor{red}
%\usepackage{ulem}
%%\allowdisplaybreaks
%\usepackage{xcolor}
%\usepackage[natbibapa]{apacite}
%\usepackage{comment}
%----------------------------------------------------------------------------------------
%	ARTICLE INFORMATION
%----------------------------------------------------------------------------------------
%\author[LIO]{Omar Kemmar}
%\author[LIO]{Karim Bouamrane}
\author[LGI2A]{Shahin Gelareh \corref{cor1}}

\cortext[cor1]{Corresponding author, shahin.gelareh@\{univ-artois.fr,gmail.com\}}
%\address[SOTON]{Centre for Operational Research,  Management Science and Information Systems (CORMSIS), University of Southampton, Southampton, SO17 1BJ, United Kingdom}
%\address[OSM]{ Portsmouth Business School, Portsmouth, United Kingdom}

% %\address[LGI2A]{Laboratoire de G\'enie Informatique et d'Automatique, Universit\'{e} Artois, F-62400, B\'{e}thune, France}

%\address[LGI2A]{LGI2A (EA 3926), Universit\'{e} d'Artois, F-62400 B\'{e}thune, France}
%\address[LIO]{Laboratoire d'informatique d'Oran, Université Oran1, BP 1524 EL Mnaouer Oran, Algerie}
\address[LGI2A]{D\'epartement R\'eseaux et T\'el\'ecommunications, Universit\'{e} d'Artois, F-62400 B\'{e}thune, France}

 %\address[LGI2A]{Laboratoire de G\'enie Informatique et d'Automatique, Universit\'{e} Artois, F-62400, B\'{e}thune, France}
% \address[NORD]{Univ Lille Nord de France, F-59000 Lille, France}
 %\address[ISIMA]{ISIMA, Universit\'{e} Blaise-Pascal, BP 10125, F-63173 Aubi\`{e}re Cedex, France}
 %\address[UBP]{LIMOS, UMR 6158-CNRS Universit\`{e} Blaise-Pascal, BP 10125, F-63173 Aubi\`{e}re Cedex, France}
% \address[IFFSTAR]{Ifsttar, Univ. Lille Nord de France, rue E\'lise\'e Reclus 20, 59666 Villeneuve d'Ascq, France}
% \address[CRIStAL]{CRIStAL, UMR 9189-CNRS, Ecole Centrale de Lille, 59651 Villeneuve d'Ascq, France}

 %\address[ITWM]{Fraunhofer Institute for Industrial Mathematics (ITWM), Kaiserslautern, Germany}

\title{A note on 'Collaborative hub location problem under cost uncertainty'}
%----------------------------------------------------------------------------------------
%	ABSTRACT
%----------------------------------------------------------------------------------------

\begin{abstract}
Three models were presented in  \emph{M.K. Khakim Habibi, Hamid Allaoui, Gilles Goncalves, Collaborative hub location problem under cost uncertainty, Computers \& Industrial Engineering Volume 124, October 2018, Pages 393-410} as models for collaborative Capacitated Multiple Allocation Hub Location Problem. In this note, we point out a few flaws in modeling. In particular, we elaborate and explain that none of the those models incorporates  any element of a collaborative activity. 
\end{abstract}
\end{frontmatter}
%----------------------------------------------------------------------------------------

%----------------------------------------------------------------------------------------
%	ARTICLE CONTENTS
%----------------------------------------------------------------------------------------
%
\section*{Introduction} % The \section*{} command stops section numbering

In the relevant literature, the term '\emph{collaboration}' is synonym for '\emph{the process of two or more people or organizations working together to complete a task or achieve a goal} (see,  Marinez-Moyano, I. J. \emph{Exploring the Dynamics of Collaboration in Interorganizational Settings}, Ch. 4, p. 83, in Schuman (Editor). Creating a Culture of Collaboration. Jossey-bass, 2006. ISBN 0-7879-8116-8.).\\

%Often in the literature, a \emph{collaborative service} is referred to a service that supports cooperative work among organizations by provision of shared access to some common resources. Therefore, there is an inherent notion of \emph{interaction} among players in this setting.\\

Three models were presented in  \cite{HABIBI2018393} as models for collaborative Capacitated Multiple Allocation Hub Location Problem.  The base model of all is one of the very early and well-studied  models of the classical hub location problem in \cite{EBERY2000614} with the same primitive assumptions which were set at the beginning of the century due to the computational limits of the CPUs in early 2000. Today, after two decades, the models of hub location have become much richer and those assumptions  are  overwritten by much more realistic features and constraints.  In this structure (see \cite{EBERY2000614}) the hub-level structure is a complete sub-graph of the resulting network and a capacity is set on the volume of flow entering a hub node. No origin-destination path travels more than one hub edge while spoke nodes are allocated to as many hub node as they wish.

\section{Literature review}

'\emph{We are particularly interested in this problem considering the presence of uncertainty due to its applicability as  well as the \underline{lack of research } literature.}\cite{HABIBI2018393}'.\\

It is often not very wise to talk about '\emph{lack of research}' in some classical topic such as location, routing and scheduling:  Interested readers are also referred to the recent HLPs literature review in \citet{Alumur&Kara2008}, \cite{campbell2012twenty}  and \citet{Farahani:2013} as well as other contributions by \citet{Campbell&Ernst&Krishnamoorthy2002} and
\citet{KaraTaner:20111}.\\

The Hub-and-spoke structures is a very active research area and lots of development. Two of the most recent and very relevant work that subsume the 'collaborative' modeling includes   \cite{GROOTHEDDE2005567} for a work on collaborative, intermodal hub networks (a case study in the fast moving consumer goods market) and \cite{coopetitive} for a coopetitive hub location model that includes both competition collaboration/cooperation. One may conclude that we are far from a  \emph{research desert} when it comes to variants of hub location problem and hub-and-spoke structures.
 
\section{Modeling and Problem Settings}

In \cite{HABIBI2018393}, the authors assume that '\emph{each supply chain is represented as a distribution network containing hub and spoke nodes}'.  The authors do not specify whether every Supply Chain $SC$ has its own origin-destination product/service to fulfill demand nodes or the same demand can be fulfilled by any of the $SC$s.\\

In \cite{HABIBI2018393}, the original model in \cite{EBERY2000614} has been referred to as 'No Collaboration' model wherein every supply chain operator operates its optimal network. Let $\chi, \alpha, \delta$ be the discount factors for collection, transfer and distribution. $W_{ij}$ represent the demand matrix and $\mathcal{N}$ stands as for the set of nodes. $C_{ij}$ and $F_k$ are the transportation cost per unit of flow and hub node setup costs, respectively. $\sigma_k^s$ is the so called supplementary cost for installing a shared hub node in scenario $s$. $(\tilde F)^s_k$ is the setup cost comprising supplementary cost $(\tilde F)^s_k=F_k+\sigma^s_k$ and $\Gamma_k$ is the capacity of hub node $k$. The variables are the following: $H_k$ is 1 if a hub is located at $k$, 0 otherwise. $Z_{ik}$ is the flow from origin $i$ to hub $k$; $Y_{kl}^i$ represents the flow from origin $i$ via hub link $k$ and $l$ and $X_{lj}^i$ represent the flow from origin $i$ to the destination $j$ via hub $l$.\\

\begin{note}\label{note:1}
It is unclear why in the 'No Collaboration', the constraints in model are indexed for a network composed of all nodes in $\mathcal{N}$. One may conclude that all the SCs are operating on the same set of nodes, i.e. $\mathcal{N}$. \\
\end{note}
The mathematical model in \cite{HABIBI2018393} follows:
\begin{align}
  {minmax}_{\forall s \in \mathcal{S}} ~~& L_s = \sum_{k\in\mathcal{N}}{F}_k H_k + \sum_{i\in \mathcal{N}}\left( \chi \sum_{k\in \mathcal{N}} C_{ik}Z_{ik} + \alpha \sum_{k\in \mathcal{N}}\sum_{l\in \mathcal{N}}C_{kl}Y_{kl}^i + \delta \sum_{l\in \mathcal{N}}\sum_{j\in \mathcal{N}}C_{ij}X_{lj}^i\right) \label{obj}\\
  s.t.: &\nonumber\\
  & \sum_{k\in \mathcal{N}} Z_{ik} = \sum_{j\in \mathcal{N}}W_{ij}, & \forall l\in \mathcal{N} \label{eq2}\\
  & \sum_{l\in \mathcal{N}} X_{lj}^i = W_{ij}, & \forall i,j\in \mathcal{N} \label{eq3}\\
  & \sum_{i\in \mathcal{N}} Z_{ik} \leq \Gamma_k H_k, & \forall k\in \mathcal{N} \label{eq4}\\
  & \sum_{l\in \mathcal{N}} Y_{kl}^i +  \sum_{j\in \mathcal{N}} X_{kj}^i = \sum_{l\in \mathcal{N}} Y_{lk}^i +  Z_{ik} & \forall i,k\in \mathcal{N} \label{eq5}\\
  & Z_{ik} \leq \sum_{j\in \mathcal{N}} W_{ij} H_k & \forall i,k\in \mathcal{N} \label{eq6}\\
  & \sum_{i\in \mathcal{N}} X_{lj}^i \leq \sum_{i\in \mathcal{N}} W_{ij} H_l & \forall l,j\in \mathcal{N} \label{eq7}\\
  & X_{lj}^i, Y_{kl}^i, Z_{ik}\geq 0 & \forall i,j,k,l\in \mathcal{N} \label{eq8}\\
  & H_k\in \{0,1\} & \forall k\in \mathcal{N} \label{eq9}
\end{align}

 In the following we review the models and elaborate on the existing flaws in modeling collaboration:

\section{Centralized Collaboration (CC)}

Again, in  the Centralized Collaboration model, the authors use the same set $\mathcal{N}$ to refer to the set of nodes for every $SC$. One may again, as in Note \autoref{note:1}, conclude that both operators are operating on exactly the same set of nodes. The flaw in this model is the following.\\

%The constraints (3) confirm that no matter which operator, it is always $w_ij$ is the supply of node $i$ for the demand of node $j$. It means that no matter what product or service, the origin-destination demand for the operators are the same. For instance if the first operator is supplying potato, and second one is supplying shoes, then the between $i$ and $j$, $w_{ij}$ ton of potato and $w_{ij}$ gigabyte of internet traffic are transported and these two product according to the constraints (2)-(9) are traveling the same origin-destination path along the network.

The objective function is the following:
\begin{align}
  {minmax}_{\forall s \in \mathcal{S}} & L_s = \sum_{k\in\mathcal{N}}\tilde{F}^s_k H_k + + \sum_{i\in \mathcal{N}}\left( \chi \sum_{k\in \mathcal{N}} C_{ik}Z_{ik} + \alpha \sum_{k\in \mathcal{N}}\sum_{l\in \mathcal{N}}C_{kl}Y_{kl}^i + \delta \sum_{l\in \mathcal{N}}\sum_{j\in \mathcal{N}}C_{ij}X_{lj}^i\right) \label{obj:CC}\\
  s.t& \nonumber\\
  & \eqref{eq2}-\eqref{eq9}\nonumber
\end{align}
subject to the same set of constraints as in \cite{EBERY2000614}, i.e. 'No Collaboration' model. In this slightly modified model, the setup cost of hub node takes form of random variable with unknown distribution but discrete probability values. \\

The authors claim that this is a model for '\emph{\underline{two networks} forming a \underline{coalition} in order to achieve \underline{shared optimal decision}}'.\\

Insofar as  the model \eqref{obj:CC}, \eqref{eq2}-\eqref{eq9}, the network node set remains the same, the origin-destination  (O-D) flows are the same as in NC and the amount of flow arriving to a node or leaving a node remains equal to the demand of every individual $SC$, which apparently operates on the entire network. \\

The least one expected to convince oneself of having two network here would be that the flow originating from a node that is a member of more than one $SC$ be the sum of supply to all those nodes (which is not the case in this model), any sign, index or constraint differentiating between the variables or constraints of every one of those 'two' networks.  \\

Nothing can be seen in this model that connect it to any other supply chain operator. This model is only in relation to itself. We are facing a single operator that is dealing with its own  uncertainty in its setup costs. This operator has not even expanded its operation over a network of another operator as the set $\mathcal{N}$ remains unchanged. Therefore, this operator works on its own network and now has to deal with its own uncertainty in setup cost. This is a unilateral decision making and has nothing to do with any sort of collaboration.\\

%The objective function \eqref{obj} and constraints \eqref{eq2}-\eqref{eq9} only deal with one isola

There are no two networks; there is no coalition and there is no share feature in anything of this mode.\\

Therefore, unless proven otherwise, this model does not represent any form or shape of collaborative  problem. As the literature is already aware of several example of similar work in CMAHLP with uncertainty in some parameters, this model by itself can hardly show any sign of novelty to be considered as a contribution.\\

\section{Centralized Collaboration with Uncertain Supplementary Cost (CCU)}
%\section{Optimized Collaboration with Uncertain Supplementary Cost(OCU)}

In \cite{HABIBI2018393}, when the minimax (CC) model  is linearized in the form of a max-regret problem, the authors refer to it as a CCU model. \\
The CCU model follows:

\begin{align}
  Min & R \label{eq11}\\
  &s.t& \nonumber\\
  & \eqref{eq2}-\eqref{eq9} \nonumber\\
  & R_s  = \sum_{k\in\mathcal{N}}\tilde{F}^s_k H_k + + \sum_{i\in \mathcal{N}}\left( \chi \sum_{k\in \mathcal{N}} C_{ik}Z_{ik} + \alpha \sum_{k\in \mathcal{N}}\sum_{l\in \mathcal{N}}C_{kl}Y_{kl}^i + \delta \sum_{l\in \mathcal{N}}\sum_{j\in \mathcal{N}}C_{ij}X_{lj}^i\right) - L_s^*   & \forall s\in\mathcal{S}\label{eq12}\\
  & R\geq R_s,  &\forall s\in S\label{eq13}\\
  &R_s,R\in \mathbb{R}\forall s\in S\label{eq14}
\end{align}

Again, beyond some algebra for getting to CCU from CC, there is no minimal link (even exchange of information on cost or any other thing) going beyond planning for anything more than one single entity (supply chain operator). \\

In all these development  a second (or any other) supply chain players is absolutely absent in any  modeling. The O-D flow are the same, the $SC$ network is the same  and we are not aware of any other $SC$ and its market share. Therefore, unless proven otherwise, this model does not represent any form or shape of collaborative problem (not even talking about location, routing or scheduling ...).

The network remains the same (i.e. $\mathcal{N}$) and this model has nothing to exchange or share with outside world. Therefore, there is no collaborative aspect in this modeling.\\

 Converting a minimax problem to a max-regret problem is closer to a textbook exercise. The literature of minmax regret model even in location problem is by far more advanced  and this model by itself can hardly be considered as a contribution given the rich body of state-of-the-art research.

\section{Optimized Collaboration with Uncertain Supplementary Cost (OCU)}
Finally, we reach to a point where the authors talk about 'other' and develop at \emph{least some constraints} for distinct pairs of $SC$s. Two additional set of variables are introduced by the authors: $I_k=1$, if a hub located at $k$ is a \emph{collaborative} one, 0 otherwise and $T_k=1$, if a hub located at $k$ is a \emph{non-collaborative} one. The model becomes:

\begin{align}
  Min & R \label{eq11}\\
  &s.t& \nonumber\\
  & \eqref{eq2}-\eqref{eq9}, \eqref{eq12}-\eqref{eq14} \nonumber\\
  &H_k = I_k+T_k   & \forall k\in\mathcal{N}\label{eq15}\\
  &Z_{ik} \leq M(1-T_k) & \forall i\in SC_a, k\in SC_b, SC_a\neq SC_b:\forall SC_a,SC_b \in \mathcal{SC} \label{eq16}\\
  &X_{lj}^i \leq M(1-T_l) & \forall i,j\in SC_a, l\in SC_b, SC_a\neq SC_b:\forall SC_a,SC_b \in \mathcal{SC} \label{eq17}\\
  & Y_{kl}^i \leq M(1-T_l) & \forall i,k\in SC_a, l\in SC_b, SC_a\neq SC_b:\forall SC_a,SC_b \in \mathcal{SC} \label{eq18}\\
  & Y_{kl}^i \leq M(1-T_k) & \forall i,l\in SC_a, k\in SC_b, SC_a\neq SC_b:\forall SC_a,SC_b \in \mathcal{SC} \label{eq19}\\
  & Y_{kl}^i \leq M(1-T_k)(1-T_l) & \forall i\in SC_a, k,l\in SC_b, SC_a\neq SC_b:\forall SC_a,SC_b \in \mathcal{SC} \label{eq20}\\
  & R_s  =  \sum_{k\in \mathcal{N}}\left( F_kT_k +\tilde{F}^s_k H_k  \right. \nonumber\\
  & \left.\sum_{k\in\mathcal{N}}(\chi \sum_{k\in \mathcal{N}} C_{ik}Z_{ik} + \alpha \sum_{k\in \mathcal{N}}\sum_{l\in \mathcal{N}}C_{kl}Y_{kl}^i + \delta \sum_{l\in \mathcal{N}}\sum_{j\in \mathcal{N}}C_{ij}X_{lj}^i)\right) - L_s^*   & \forall s\in\mathcal{S}\label{eq21}\\
  & I_k, T_k \in \{0,1\}, \forall k\ \in \mathcal{N}\label{eq22}
\end{align}
\begin{note}
According to the constraints  $SC_a\subset \mathcal{N}$ and $SC_b\subset \mathcal{N}$, $SC_a \neq SC_b$ for every $SC_a, SCb \in SC$ and we assume it is implicitly said that $\bigcup_{a\in SC}SC_a=\mathcal{N}$ and $SC_a\cap SC_b, \forall ~~a,b$ may be non-empty. This contradicts in what has been concluded in Note \autoref{note:1} because over there, every $SC$ for example $SC_a$ was solving the model on the entire $\mathcal{N}$ and not only on nodes of $SC_a$.\\
\end{note}

But given the new notations let us suppose that the node set $\mathcal{N}$ is composed of union of node sets for every $SC_a \in SC$ nodes. This means that the authors basically have made up some subsets of $\mathcal{N}$ with possibly non-empty intersections (i.e. the sets are not necessarily mutually exclusive) and imposed additional constraints on the CCU model to achieve OCU.

There are a few flaws which are elaborated in the following:
\subsection{The objective value of OCU}

\begin{thm}\label{thm1}
The objective function of OCU can never be any better than the objective of CCU.
\end{thm}
\begin{pf}
In simple words and without resorting to any complicated polyhedral theory and descriptions, one notices that OCU is no less constrained  than model CCU and therefore the number of feasible solutions in OCU (in the space of $X,Y,Z$) is no more than those in CCU. In particular, given that $SC_a \neq SC_b$ for every $SC_a, SC_b \in SC$, one may say that  OCU is in fact strictly more constrained  than model CCU. As such, the set of feasible solutions in OCU is strictly contained in the set of feasible solutions of CCU (in the space of $X,Y,Z$) and no objective function delivers a better objective value on $\mathcal{P}(OCU)$ than on $\mathcal{P}(CCU)$, being polytopes of the corresponding problems. As such, any claim and computational results claiming otherwise is provably wrong.
%Let $\mathcal{P}(OCU)$ and $\mathcal{P}(CCU)$ be the polytopes associated with the OCU and CCU models. Clearly, $\mathcal{P}(OCU)\subseteq \mathcal{P}(CCU)$. More precisely, as $SC_a \neq SC_b$ for every $SC_a, SCb \in SC$, $\mathcal{P}(OCU)\subset \mathcal{P}(CCU)$. As such, no objective function achives a better objective value on \mathcal{P}(OCU) than on \mathcal{P}(CCU). Therefore any claim
\end{pf}

Therefore, the so-called '\emph{collaboration}', if any, is an absolute loss for the supply chain as a whole. The authors would not even need to conduct computational experiments to realize this fact. Moreover, any numerical outcome that contradicts this, is most probably a result of mistake and error, in an optimistic viewpoint.\\

\subsection{flaw in modeling:\#1}
According to the constraints \eqref{eq16}-\eqref{eq19}, no two subsets $SC_a$ and $SC_b$ are exactly the same. It is unclear whether the intersection is non-empty or not. Suppose $k$ is a node in the intersection of $SC_a$ and $SC_b$ and is a hub node. The first question arises here, is the following: who has paid the setup cost for such a hub node, was it $SC_a$ or $SC_b$? Which stakeholder in this problem description, if any, is in charge of setting up the hub nodes? How can one say that this hub belongs to $SC_a$ but not to $SC_b$, or the other way around? \\

In addition, according to the constraints \eqref{eq16}-\eqref{eq19}, none of these two $SC$s are authorized to use this hub node unless $T$ variable becomes 0 (what authors call it a variable for determining 'non-collaborative' nature of a hub node).  First of all, for such a node in the intersection of node sets, we do not know to whom this hub belongs. Even more, if that node is also in the node set of any other $SC$, the later cannot use it either, unless pay extra to turn $I_k$ to 1.
\\

If we assume that the intersection of every pair of $SC_a$ and $SC_b$  is an empty set, the issue becomes even more serious. Suppose that $SC_a$ and $SC_b$ operating on disjoint node sets and they do not have any node in common. In this case, how can one say that there has been a flow from $i\in SC_a$ to $j\in SC_b$, i.e. $W_{ij}\neq 0$. How this flow  even existed when every operator was operating on its own node set? Where does this ghost O-D flow come from.\\

To that one can add the following: Given constraints \eqref{eq18}-\eqref{eq20}, this model does not deliver any optimal solution in, which for a hub node $k$, $T_k=1$.

In any case, it is clear that some unknown and absent stakeholder, as an autocracy, is present who forcefully wants all these nodes to be regrouped in one decision and one network, magically generates additional market share among the pairs that did not know each others before and solve one problem as a whole.  This ghost stakeholder has only one goal, he does not want to see many players, just wants one single network.  \\

Given that, there is absolutely no collaboration going on here.

\subsection{flaw in modeling:\#2}
An average  reader understands that the OCU model is solved once and gives solution for the whole merger, otherwise it does not make sense that every $SC$ solve one separate OCU and pays the full cost of setup for a hub node that is in the intersection of two $SC$s and is already paid by one or more other $SC$. But if so, why we do not see any additional flow? The same $W$ matrix as before is used in OCU. The volume in this merger has not increased by summing up the volume of all $SC$s and the capacities remain the same even in this merger!. A single $SC$ given its volume of flow, set $\Gamma_k$ capacity on hub $k$ and when they all merge and volume must normally increase nothing happens to the volumes and structures.\\

\subsection{flaw in modeling:\#3}
The nonlinear constraints \eqref{eq20} are totally irrelevant. They are already implied by the constrains \eqref{eq18}-\eqref{eq19}.  Therefore, all trivial linearization algebra that follows, including the 'Proposition 1' in \cite{HABIBI2018393} are irrelevant and over engineered.

\subsection{flaw in modeling:\#4}
The definition of variables $T_{k}$ and $I_k$ is curious:   

\begin{align}
  I_k \left\{
           \begin{array}{ccc}
             1 & \mbox{if a hub located at $k$ is collaborative} \\
             0 & \mbox{oterwise}. \\
           \end{array}
         \right.\hspace{2cm}
         T_k \left\{
           \begin{array}{ccc}
             1 & \mbox{if a hub located at $k$ is non-collaborative} \\
             0 & \mbox{oterwise}. \\
           \end{array}
         \right.
\end{align}

Given constraints \eqref{eq15}, i.e. $H_k=I_k+T_k$, if any of $I_k$ or $T_k$ takes 1, it means that $k$ is a hub node. Anyone of $I$ or $T$ would implies the other in the definition, so why should one declare two sets of variables and unnecessarily enlarge the polyhedral description in a higher dimension? In particular, just in the context of modeling (not approving the model itself)  definition of variable $I_k$ apparently serves no purpose and $I_k$ is already implied by $T_k$.

%It is amazing that before even  any network merging takes place in the form of $\bigcup_{a\in SC}SC_a=\mathcal{N}$, there is demand for supply of node $i\in SC_a$ and $j\in SC_b$! As it is explicitly said in constraints \eqref{eq16}-\eqref{eq19}, no two subsets $SC_a$ and $SC_b$ are exactly the same. How comes that a node in $SC_a$, namely $i$ was aware of its supply to node $j$ in $SC_b$? $SC$s were supposed to be independent operators operating their own networks before they merge, as far as a average reader would understand.\\

%\subsection{flaw in modeling:\#4}
%Assume that $i\in SC_a$ is the source of flow $W_{ij}$ and $j\in SC_b$ is the closest node to $i$. According to the constraints \eqref{eq16}, $i$ could not send its flow to $j$ unless  $j$ becomes a so-called 'collaborative' hub node otherwise it has to reroute it in a significantly longer path with higher cost for the pleasure of '\emph{collaboration}'. It is actually possible, if this is compensated by a feature in modeling representing collaboration.

\subsection{flaw in modeling:\#5}
Cost sharing is not a post-processing phase. The model and solution must show that there is an added value in the collaboration. We know that OCU never generates any better feasible solution than any other model mentioned. The question is, when no positive sign is showing up, it means that the overall merger is loosing, while $SC$s are loosing at different levels. \\

It is curious that when the overall cost increasing, no additional market is being generated as a result of the merger, no economy of scale is exploited by deploying more efficient (perhaps larger with higher capacity) transporter on any hub-level link, or anything else, this 'merger' can make any sense for anyone or would mean anything similar to a '\emph{collaboration}'. \\

The optimization must be made in presence of a motivating factor such as reduction of overall cost due to increase in market share, reduce in transportation cost due to the establishment of new hub-to-hub links or something that must become better through the so called '\emph{collaboration}'. 
%
%\subsection{Stakeholder issue}
%Participants in any collaboration are autonomous independent entities which decide on the

\section{Conclusion}
None of the proposed models have much to do with any sort of 'collaboration'. It is very trivial, as shown above, that the heart of this article which is the OCU model does not deliver any better feasible solution. The proposed model can hardly have any real-life realisation. The proposed model OCU increases the overall cost of operation (perhaps gain for some but loose for the overall merger) and the postprocessing phase, rather than being a cost sharing mechanism is a disaster sharing technique. The model is basically promoting the notion of '\emph{loosing together and crying together}'.

\bibliographystyle{apalike}
\bibliography{sample}

\begin{thebibliography}{}

\bibitem[Alumur and Kara, 2008]{Alumur&Kara2008}
Alumur, S. and Kara, B.~Y. (2008).
\newblock Network hub location problems: the state of the art.
\newblock {\em European Journal of Operational Research}, 190:1--21.

\bibitem[Campbell et~al., 2002]{Campbell&Ernst&Krishnamoorthy2002}
Campbell, J.~F., Ernst, A.~T., and Krishnamoorthy, M. (2002).
\newblock Hub location problems.
\newblock In Drezner, Z. and Hamacher, H.~W., editors, {\em Facility Location:
  {A}pplications and {T}heory}, pages 373--407. Springer.

\bibitem[Campbell and O'Kelly, 2012]{campbell2012twenty}
Campbell, J.~F. and O'Kelly, M.~E. (2012).
\newblock Twenty-five years of hub location research.
\newblock {\em Transportation Science}, 46(2):153--169.

\bibitem[Ebery et~al., 2000]{EBERY2000614}
Ebery, J., Krishnamoorthy, M., Ernst, A., and Boland, N. (2000).
\newblock The capacitated multiple allocation hub location problem:
  Formulations and algorithms.
\newblock {\em European Journal of Operational Research}, 120(3):614 -- 631.

\bibitem[Farahani et~al., 2013]{Farahani:2013}
Farahani, R.~Z., Hekmatfar, M., Arabani, A.~B., and Nikbakhsh, E. (2013).
\newblock Survey: Hub location problems: A review of models, classification,
  solution techniques, and applications.
\newblock {\em Comput. Ind. Eng.}, 64(4):1096--1109.

\bibitem[Groothedde et~al., 2005]{GROOTHEDDE2005567}
Groothedde, B., Ruijgrok, C., and Tavasszy, L. (2005).
\newblock Towards collaborative, intermodal hub networks: A case study in the
  fast moving consumer goods market.
\newblock {\em Transportation Research Part E: Logistics and Transportation
  Review}, 41(6):567 -- 583.
\newblock Global Logistics.

\bibitem[Habibi et~al., 2018]{HABIBI2018393}
Habibi, M.~K., Allaoui, H., and Goncalves, G. (2018).
\newblock Collaborative hub location problem under cost uncertainty.
\newblock {\em Computers \& Industrial Engineering}, 124:393 -- 410.

\bibitem[Kara and Taner, 2011]{KaraTaner:20111}
Kara, B.~Y. and Taner, M.~R. (2011).
\newblock Hub location problems: {T}he location of interacting facilities.
\newblock In Eiselt, H.~A. and Marianov, V., editors, {\em Foundations of
  location analysis}, pages 273--288. Springer.

\bibitem[Monemi et~al., 2017]{coopetitive}
Monemi, R., Gelareh, S., Hanafi, S., and Maculan, N. (2017).
\newblock A co-opetitive framework for the hub location problems in
  transportation networks.
\newblock {\em Optimization}.

\end{thebibliography}

%----------------------------------------------------------------------------------------

\end{document}